\documentclass[a4paper,12pt]{article}

\usepackage{amsmath}
\usepackage{amssymb}
\usepackage{amsfonts}
\usepackage{amsthm}
\usepackage{mathrsfs}
\usepackage{graphicx}
\usepackage{makeidx}
\usepackage{subcaption}
\usepackage{color}
\usepackage{mdframed}
\usepackage{ulem}

\usepackage[left=3cm,right=3cm,top=3cm,bottom=3cm]{geometry}

\newcommand{\R}{\mathbb{R}} 
\newcommand{\N}{\mathbb{N}} 

\newtheorem{theorem}{Theorem}
\newtheorem{lemma}{Lemma}

\newtheorem{definition}{Definition}
\newtheorem{assumption}{Assumption}
\newtheorem{remark}{Remark}

\newcommand{\new}[1]{{\color{black}#1}}

\newcommand{\bg}{{\bf g}}
\newcommand{\bx}{{\bf x}}
\newcommand{\by}{{\bf y}}

\newcommand{\sF}{{\mathscr F}}
\newcommand{\sM}{{\mathscr M}}
\newcommand{\sN}{{\mathscr N}}
\newcommand{\sP}{{\mathscr P}}
\newcommand{\sQ}{{\mathscr Q}}
\newcommand{\sR}{{\mathscr R}}
\newcommand{\sS}{{\mathscr S}}
\newcommand{\sX}{{\mathscr X}}
\newcommand{\sY}{{\mathscr Y}}
\newcommand{\mom}{\text{mom}}
\newcommand{\sos}{\text{sos}}

\title{\bf Geometry of exactness of moment-SOS relaxations for polynomial optimization}

\begin{document}

\author{Didier Henrion$^{1,2}$}

\footnotetext[1]{CNRS; LAAS; Universit\'e de Toulouse, 7 avenue du colonel Roche, F-31400 Toulouse, France. }
\footnotetext[2]{Faculty of Electrical Engineering, Czech Technical University in Prague,
Technick\'a 2, CZ-16626 Prague, Czechia.}

\date{Draft of \today}

\maketitle

\begin{abstract}
The moment-SOS (sum of squares) hierarchy is a powerful approach for solving {approximately and} globally non-convex polynomial optimization problems (POPs) at the price of solving a family of convex semidefinite optimization problems (called moment-SOS relaxations) of increasing size, controlled by an integer, the relaxation order. We say that a relaxation of a given order is exact if solving the relaxation actually solves the POP globally. In this note, we study the geometry of the exactness cone, defined as the set of polynomial objective functions for which the relaxation is exact. Generalizing previous foundational work on quadratic optimization on real varieties, we prove by elementary arguments that the exactness cones are unions of semidefinite representable cones monotonically embedded for increasing relaxation order.
\end{abstract}

\section{Solving POPs with the moment-SOS hierarchy}

Given a compact semialgebraic set $\sX  {\:\subset\:} \R^n$ and a polynomial $f$ in the vector space $\R[\bx]_d$ of polynomials {in} $\bx  \in \R^n$ of degree up to $d$, consider the \emph{polynomial optimization problem} (POP)
\begin{equation}\label{pop}
v(f):= \min_{\bx \in \sX} f(\bx).
\end{equation}
 The notation  emphasizes that the optimal value depends parametrically on the objective function.

The key observation behind the moment-SOS (sum of squares) hierarchy \cite{l15,hkl20,h23,n23} is that the POP is equivalent to the primal-dual problems
\[
 \begin{array}{rcllcll}
 v(f) & = & \min_{\by} & \ell_{\by}(f) & = & \max_{v} & v\\
&& \mathrm{s.t.} & \by \in \sM(\sX)_d & & \mathrm{s.t.} & f-v \in \sP(\sX)_d \\
&&& \ell_{\by}(1) = 1
\end{array}
\]
where the dual maximization problem consists of finding the largest lower bound $v \in \R$ on $f$ on $\sX$, formalized as a linear conic problem in the convex cone $\sP(\sX)_d$ of polynomials of degree up to $d$ which are positive on $\sX$. The primal minimization problem is over vectors $\by$ in the convex cone $\sM(\sX)_d$ of moments of degree up to $d$ of positive measures on $\sX$, which is the convex dual (i.e. the set of positive linear functionals) of $\sP(\sX)_d$, cf. e.g. \cite[Thm. 8.1.2]{n23}. The linear objective function in the primal conic problem is $\ell_{\by}(f) := \int_{\sX} f(\bx)d\mu(\bx)$ where $\mu$ is a positive measure with moment vector $\by$, and the linear constraint $\ell_{\by}(1)=\int_{\sX} d\mu(\bx)=1$ enforces that $\mu$ is a probability measure.

If $\sX:=\{\bx \in \R^n : g_i(\bx)\geq 0, i=1,\ldots,m\}$ is basic semialgebraic, defined by a given polynomial vector $\bg=(g_i)_{i=1,\ldots,m}$, then $\sP(\sX)_d$ can be approximated {from inside} with other convex cones, the truncated quadratic modules
\[
\sQ(\bg)^r_d:=\{p \in \R[\bx]_d : p = \sum_{i=0}^{m} s_i g_i, \: s_i \in \Sigma[\bx], \: s_i g_i \in \R[\bx]_{2r}\}
\]
where $\Sigma[\bx] \subset \R[\bx]$ is the convex cone of sums of squares (SOS) of polynomials of $\bx$ and  $g_0(\bx):=1$.

\new{
Since set $\sX$ is bounded, it is always possible to add a redundant ball constraint to its description. So the following assumption is without loss of generality.
\begin{assumption}\label{compact}
Let $R > 0$ be the radius of an Euclidean ball including $\sX$. It holds $R^2-\sum_{k=1}^n x^2_i \in \sQ(\bg)^r_2$ for all $r \in \N$.
\end{assumption}
}

\begin{remark}\label{equal}
Note that if a polynomial equation enters the definition of $\sX$, i.e. $g_i(\bx) = 0$ instead of $g_i(\bx) \geq 0$ for some $i=1,\ldots,m$, then the corresponding weight $s_i \in \R[\bx]$ in the quadratic module $\sQ(\bg)^r_d$ is not constrained in sign, while satisfying $s_ig_i \in \R[\bx]_{2r}$. {Recalling the identity $s_i = ((s_i+1)/2)^2-((s_i-1)/2)^2$,  every polynomial of $\R[\bx]$ can be written as the difference of two polynomials of $\Sigma[\bx]$.} This is consistent with the fact that two inequalities of opposite signs are equivalent to an equation. Without loss of generality, and for notational conciseness, in this note we use only inequalities.
\end{remark}

Note that by construction the quadratic modules are monotonically embedded for decreasing relaxation order:
\begin{equation}\label{inclusion}
\sQ(\bg)^r_d \subset \sQ(\bg)^{r+1}_d \subset \sP(\sX)_d.
\end{equation}
Contrary to $\sP(\sX)_d$, the truncated quadratic module $\sQ(\bg)^r_d$ is \emph{semidefinite representable}, i.e. it is the linear projection of a spectrahedron, itself defined as a linear section of the cone of positive semidefinite quadratic forms. Practically, this means that linear optimization in $\sQ(\bg)^r_d$ can be done efficiently with powerful interior-point algorithms. 

Let us denote by $\sR(\bg)^r_d$ the convex cone dual to the truncated quadratic module $\sQ(\bg)^r_d$.
By convex duality, for the primal problem we have the reversed monotone embedding
$\sR(\bg)^r_d \supset \sR(\bg)^{r+1}_d \supset \sM(\bg)_d$ meaning that the moment cone is approximated from outside, or relaxed. This motivates the terminology \emph{moment relaxation} to refer to  $\sR(\bg)^r_d$.

Now we have all the ingredients to define the \emph{moment-SOS hierarchy} also known as Lasserre's hierarchy: a family of primal-dual convex semidefinite optimization problems whose size is controlled by the relaxation order $r \in \N$:
\begin{equation}\label{momsos}
\begin{array}{rcllcrcll}
\mom(f)^r & := & \inf_{\by} & \ell_{\by}(f) & \geq & \sos(f)^r & := & \sup_{v} & v\\
&& \mathrm{s.t.} & \by \in \sR(\bg)^r_d & & & & \mathrm{s.t.} & f-v \in \sQ(\bg)^r_d \\
&&& \ell_{\by}(1) = 1.
\end{array}
\end{equation}
This primal-dual pair of semidefinite optimization problems is called the \emph{moment-SOS relaxation} of order $r$. Note that by construction
\begin{equation}\label{bounds}
\sos(f)^r \leq \mom(f)^r \leq v(f).
\end{equation}

\new{Under Assumption \ref{compact}},
it follows from \cite{jh16} that in \eqref{momsos} the primal is attained (i.e. the infimum is a minimum), there is no duality gap (i.e. the infimum equals the supremum) and the relaxed values are monotonically converging lower bounds on the value:
\[
\sos(f)^r = \mom(f)^r \leq \sos(f)^{r+1} = \mom(f)^{r+1} \leq \sos(f)^ {\infty} = \mom(f)^{\infty} = v(f)
\]
{where $\sos(f)^{\infty} := \lim_{r\to\infty} \sos(f)^r$ and likewise for $\mom(f)^{\infty}$.}

\section{Exactness cone}

Beyond asymptotic convergence guarantees, it is important to know whether the moment-SOS relaxation of a given order $r$ is \emph{exact}, i.e. whether $\sos(f)^r=v(f)$. If this is the case, there is no need to increase $r$ and solve larger semidefinite optimization problems.

In this note, we are interested in the geometry of the \emph{exactness cone}, defined as the set of objective functions which are such that the moment-SOS relaxation \eqref{momsos} is exact.

\begin{definition}\label{exactnesscone}
The exactness cone of degree $d$ at relaxation order $r$ is defined by
\[
\new{\sF(\bg)^r_d:=\{f \in \R[\bx]_d : f-v(f) \in \sQ(\bg)^r_d\}}.
\]
\end{definition}

Note that this set is a cone since $v(af)=av(f)$ for all $a \geq 0$.

Our main result states that the exactness cone is a (generally uncountable and non-convex) union of semidefinite representable cones. We also describe the convex geometry of the exactness cones, their monotone embedding, and how they are related to normal cones of the moment relaxations. 

Our analysis is elementary. It is inspired by the foundational work \cite{chs21} which focused on the particular case of the Shor relaxation ($r=1$) with $f$ linear ($d=1$) or quadratic ($d=2$), and $\sX$ a real algebraic variety defined by quadratic equations. We believe that our contribution consists of considerably simplifying and extending {the geometric description of the exactness cones given in \cite{chs21}} to higher order relaxations of general semialgebraic sets.

\section{Main result}

\begin{theorem}\label{thm}
	\new{Under Assumption \ref{compact},}
the exactness cone at relaxation order $r$ is given by
\begin{equation}\label{exact}
\sF(\bg)^r_d = \bigcup_{\hat{\bx}\in \sX} \sS_{\hat{\bx}}(\bg)^r_d
\end{equation}
where each cone
\[
\sS_{\hat{\bx}}(\bg)^r_d := \{f \in \R[\bx]_d : f-f(\hat{\bx}) \in \sQ(\bg)^r_d\}
\]
is semidefinite representable.
\end{theorem}

{\bf Proof:}
If $f \in \sF(\bg)^r_d$ then $\sos(f)^r=v(f)=f(\hat{\bx})$ where $\hat{\bx}$ is a global minimizer of $f$ on $\sX$. As $\sos(f)^r$ is attained, we have $f-\sos(f)^r=f-f(\hat{\bx}) \in \sQ(\bg)^r_d$ and thus $f$ belongs to $\sS_{\hat{\bx}}(\bg)^r_d$ for some $\hat{\bx} \in \sX$.

Conversely, assume $f \in \sS_{\hat{\bx}}(\bg)^r_d$ for some $\hat{\bx} \in \sX$. Then from \eqref{bounds} we have $f(\hat{\bx}) \leq \sos(f)^r \leq \mom(f)^r \leq v(f) \leq f(\hat{\bx})$, which implies that $\sos(f)^r = \mom(f)^r = v(f)$ and hence $f \in \sF(\bg)^r_d$.

{The} cone  $\sS_{\hat{\bx}}(\bg)^r_d$ is semidefinite representable since it is the projection of linear sections of the SOS cone $\Sigma[\bx]$, which is itself semidefinite representable.
 $\Box$

\begin{remark}
The truncated quadratic module $\sQ(\bg)^r_d$ can be replaced by any other semidefinite representable approximation of $\sP(\sX)$, e.g. the preordering of $\bg$, or any other Positivstellensatz \cite{l15,n23}.
\end{remark}

\begin{remark}
The exactness cone \new{of order $r$} is semialgebraic since $f$ belongs to $\sF(\bg)^r_d$ whenever 1) $f$ belongs to $\sQ(\bg)^r_d$, a semidefinite representable hence semialgebraic cone, and 2) there exists $\hat{\bx}$ in $\sX$, a semialgebraic set, such that $f$ vanishes at $\hat{\bx}$.
\end{remark}

\begin{remark}
It follows from the proof of Theorem \ref{thm} that  it is enough to restrict the union \eqref{exact} to points $\hat{\bx} \in \sX$ which are optimal for some objective function $f$. If $d=1$, the union can be restricted to the set of extreme points of the convex hull of $\sX$, {which is however not easy to describe explicitly in general.} {If $d\geq 2$, an anonymous reviewer noticed that every point $\hat{\bx} \in \sX$ is the unique global minimizer of the quadratic objective function $f(\bx)=(\bx-\hat{\bx})^T (\bx-\hat{\bx})$, which implies that the union \eqref{exact} cannot be reduced to a subset of $\sX$.}
\end{remark}

\begin{remark}
In our definition of the exactness cone, we \new{require that the dual value is attained in \eqref{momsos}. This is called SOS exactness in \cite{b22,bm20}.} We may relax our attainment requirements. The corresponding exactness cones would be slightly larger, at the price of more technicalities.
\end{remark}

\section{Geometry of exactness cones}

\begin{lemma}\label{embed}
\new{Under Assumption \ref{compact},} the exactness cones are monotonically embedded for increasing relaxation order:
\[
\R[\bx]_0 \subset \sF(\bg)^r_d \subset \sF(\bg)^{r+1}_d \subset \overline{\sF(\bg)^{\infty}_d}=\R[\bx]_d.
\]
\end{lemma}

{\bf Proof}: 
The inclusion ${\R}[\bx]_0 \subset \sF(\bg)^r_d$ follows from the translation invariance of the value: $v(f+a)=v(f)+a$ for all $a \in \R$. The inclusion $\sF(\bg)^r_d\subset \sF(\bg)^{r+1}_d$ follows from the inclusion relations \eqref{inclusion}. Finally, the identity $\overline{\sF(\bg)^{\infty}_d}=\R[\bx]_d$ follows from Putinar's Positivstellensatz -- see e.g. \cite[Thm. 2.15]{l15} -- which states that \new{under Assumption \ref{compact}} every degree $d$ polynomial strictly positive  on $\sX$ belongs to $\sQ(\bg)^r_d$ for sufficiently large $r \in \N$, i.e. $\overline{\sQ(\bg)^{\infty}_d}=\sP(\sX)_d$ {and by duality $\overline{\sR(\bg)^{\infty}_d}=\sM(\sX)_d$. From the definition of the exactness cone and the global minimum $\hat{v}$ on $\sX$, it holds $\overline{\sF(\bg)^{\infty}_d}:=\{f \in \R[\bx]_d : \ell_{\hat{\by}}(f)=\hat{v}=v(f), \:\hat{\by} \in  \sM(\sX)_d, \:f-\hat{v} \in \sP(\sX)_d\}=\R[\bx]_d$.}
$\Box$

\begin{remark}\label{closure}
Note that in Lemma \ref{embed} the closure is required for asymptotic exactness of all degree $d$ polynomial objective functions, i.e. $\overline{\sF(\bg)^{\infty}_d}=\R[\bx]_d$. A classical example is $n=m=1$, $g_1(x)=-x^2$ for which $\sX=\{0\}$. Whereras $f(x)=\pm x \notin \sQ(\bg)^r_1$ for finite $r \in \N$, it holds that $f(x)+\varepsilon=\frac{\varepsilon}{2}+\frac{\varepsilon}{2}(1\pm\frac{x}{\varepsilon})^2-\frac{x^2}{2\varepsilon}  \in \sQ(\bg)^{\new{2}}_1$ for every $\varepsilon>0$, see \cite[Ex. 1.3.4]{b22} or \cite[\S 2.5.2]{n23}. For this example $\sF(\bg)^r_1=\R[\bx]_0$ for all finite $r$, and $\overline{\sF(\bg)^{\infty}_1}=\R[\bx]_1$.
\end{remark}

{The vector space of polynomials $\R[\bx]_d$ can be identified with $\R^{n+d \choose n}$. Using the linear functional $\ell_{\by}(f)$, it is in duality with the space of vectors $\by \in \R^{n+d \choose n}$.}
Given $\hat{\bx} \in \sX$, let $\by_{\hat{\bx}} {\:\in\: \R^{n+d \choose n}}$ be the Dirac vector at $\hat{\bx}$, i.e. such that $\ell_{\by_{\hat{\bx}}}(f)=f(\hat{\bx})$ for every $f \in \R[\bx]_d$.
Given a convex set $\sY {\:\subset\:} \R^N$, let {$\sN_{\sY}(\hat{\by}):=\{f \in \R[\bx]_d : \ell_{\by-\hat{\by}}(f) \leq 0, \: \forall \:\by \in \sY\}$ denote the normal cone to $\sY$ at a point $\hat{\by} \in \sY$.} 
With these notations, we have the following geometric counterpart to Theorem \ref{thm}.

\begin{theorem}\label{normal}
Up to the sign \new{and closure}, the spectrahedral cones of Theorem \ref{thm} are normal cones to the moment relaxation:
\[
\sN_{\sR(\bg)^r_d}(\by_{\hat{\bx}})=\new{-\overline{\sS_{\hat{\bx}}(\bg)^r_d}}
\]
for all $\hat{\bx} \in \sX$.
\end{theorem}

{\bf Proof:} 
{From the definition of the normal cone: $-\sN_{\sR(\bg)^r_d}(\by_{\hat{\bx}}) = \{ f \in \R[\bx]_d : \ell_{\by}(f) \geq \ell_{\by_{\hat{\bx}}}(f)=f(\hat{\bx}) \:\: \forall \by \in \sR(\bg)^r_d\} = \{ f \in \R[\bx]_d : \ell_{\by}(f-f(\hat{\bx})) \geq 0 \:\: \forall \by \in \sR(\bg)^r_d\}$ which is the same as $\{ f \in \R[\bx]_d : f-f(\hat{\bx}) \in \new{\overline{\sQ(\bg)^r_d }}\}
 = \new{\overline{\sS_{\hat{\bx}}(\bg)^r_d}}$ since \new{$\overline{\sQ(\bg)^r_d}$ is the dual cone to $\sR(\bg)^r_d$.}
$\Box$

{Theorem} \ref{normal} is a generalization of \cite[Prop. 4.7]{chs21} which considers only the case $d\leq2$ and the Shor relaxation (i.e. $r=1$) for real varieties $\sX$ generated by quadratic polynomials.

\new{If $\sX$ has an interior point, then the truncated quadratic module $\sQ(\bg)^r_d$ is closed, and the dual is attained in \eqref{momsos} (i.e. the supremum is a maximum). 
If $\sX$ does not have an interior point, e.g. if it is a low-dimensional algebraic variety, then additional algebraic or geometric conditions are required for the dual to be attained, cf. \cite{bm20,b22,n23}. This is why the closure is required in the statement of Theorem \ref{normal}.}

It is of interest to know whether we have exactness at relaxation order $r$ for all objective functions. \new{A sufficient} condition is that $\sR(\bg)^r_d=\sM(\sX)_d$, i.e. all vectors of the moment relaxation are convex combinations of Dirac vectors of $\sX$.
If $d=1$ the condition becomes $\sR(\bg)^r_1 = \R_+ \times \mathrm{conv}\:\sX$,
the Cartesian product of the positive real line with the convex hull of $\sX$.
%
%

\section{Examples}

{\bf Data availability statement:} The Matlab codes to generate the figures in the following examples are available upon request.

\subsection{Finite set}

Let us revisit \cite[Ex. 4.4]{chs21} where $d=1$ and
\[
\begin{array}{rcl} 
\sX  & = & \{\bx \in \R^2 : g_1(\bx)=2x_2-2x^2_2+x_1x_2=0, \: g_2(\bx)=-x_1+x_2+x^2_1-x^2_2=0\} \\
& = & \{(0,0),(0,1),(1,0),(2,2)\}
\end{array}
\]
is a finite set consisting of four points.
Then the exactness cone $\sF(\bg)^r_d$ is the union of four semidefinite representable cones. For example if $d=r=1$ one of these cones is
\[
\begin{array}{rcl}
\sS(2,2)^1_1 &  = & \{f \in \R[\bx]_1 : f(\bx)-f(2,2) = s_0(\bx) + s_1 g_1(\bx) + s_2 g_2(\bx), \\
& & \quad s_0 \in \Sigma[\bx] \cap \R[\bx]_2, s_1 \in \R, s_2 \in \R\}
\end{array}
\]
where $f(\bx)= f_0 + f_1 x_1 + f_2 x_2$,
which can be written more explicitly by expressing the quadratic SOS constraint
\[
s_0(\bx) = (1,x_1,x_2) X (1,x_1,x_2)^T, \: X = \left(x_{ij}\right)_{i,j=1,2,3} \succeq 0
\]
with a 3-by-3 positive semidefinite Gram matrix $X$, and identifying like powers of $\bx$ in the equation $f(\bx)-f(2,2) = s_0(\bx) + s_1 g_1(\bx) + s_2 g_2(\bx)$:
\[
\begin{array}{rcl}
-2f_1-2f_2 & = & x_{11} \\
f_1 & = & 2x_{21}-s_2 \\
f_2 & = & 2x_{31}+2s_1+s_2 \\
0 & = & x_{22} + s_2 \\
0 & = & 2x_{32} + s_1 \\
0 & = & x_{33} - 2 s_1 - s_2.
\end{array}
\]
Therefore
\[
\sS(2,2)^1_1 = \R \oplus \{(f_1,f_2) \in \R^2 :  \left(\begin{array}{ccc}
-4f_1-4f_2 & \new{f_1+s_2} & \new{f_2-2s_1-s_2} \\
f_1+s_2 & -2s_2 & \new{s_1} \\
f_2-2s_1-s_2 & s_1 & 4s_1+2s_2
\end{array}\right) \succeq 0, (s_1,s_2) \in \R^2\}
\]
is a projection of a 4-dimensional cubic spectrahedral cone.

\begin{figure}
	\centering
	\begin{minipage}{0.48\textwidth}
		\includegraphics[width=\textwidth]{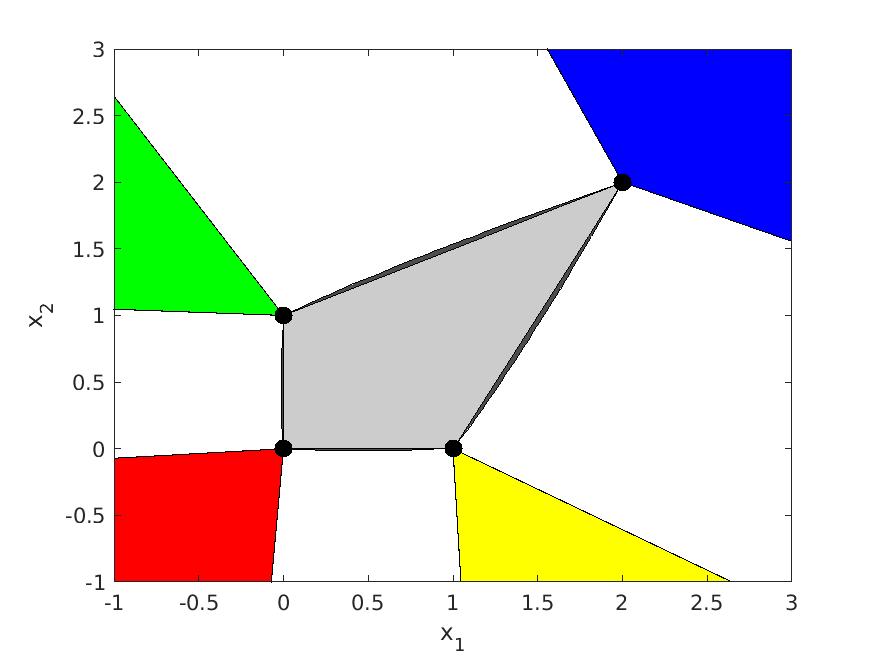}
		\subcaption{\footnotesize First moment relaxation and  normal cones}
	\end{minipage}
	\begin{minipage}{0.48\textwidth}
		\includegraphics[width=\textwidth]{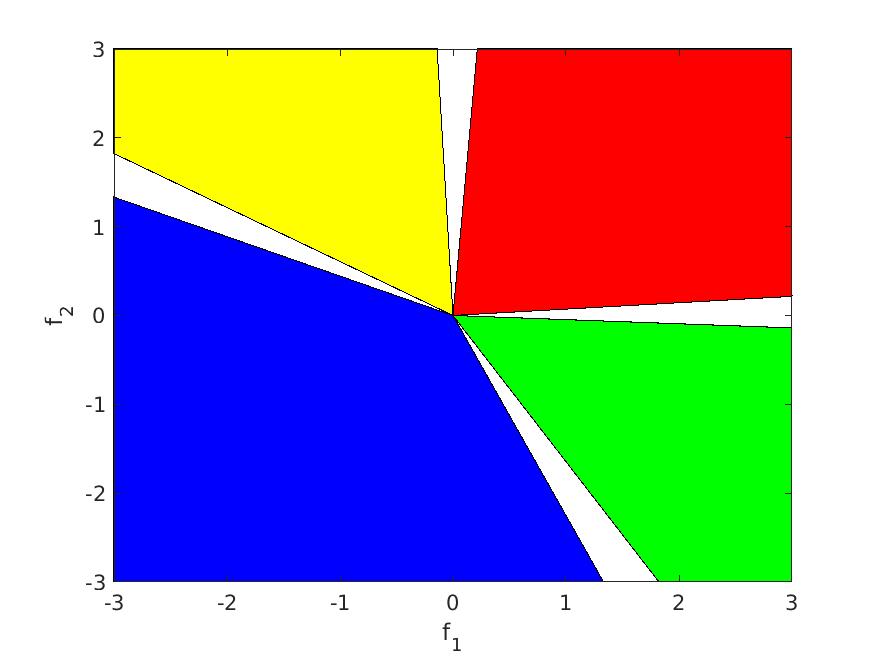}
		\subcaption{\footnotesize First exactness cone}
	\end{minipage}
	\caption{\label{fig:fourpoint1} Finite set. Left (a): first moment relaxation (dark gray) including the convex hull (light gray) of $\sX$ (four black points), and normal cones at the points (colored). Right (b): first exactness cone (colored).}
\end{figure}

\begin{figure}
	\centering
	\begin{minipage}{0.48\textwidth}
		\includegraphics[width=\textwidth]{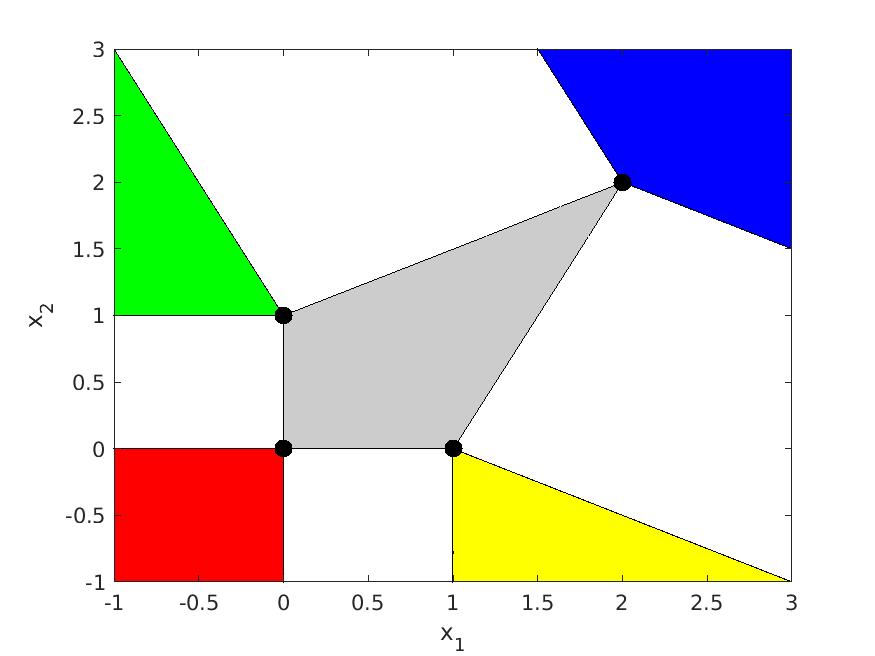}
		\subcaption{\hspace{-.21em}\footnotesize Second moment relaxation and  normal cones}
	\end{minipage}
	\begin{minipage}{0.48\textwidth}
		\includegraphics[width=\textwidth]{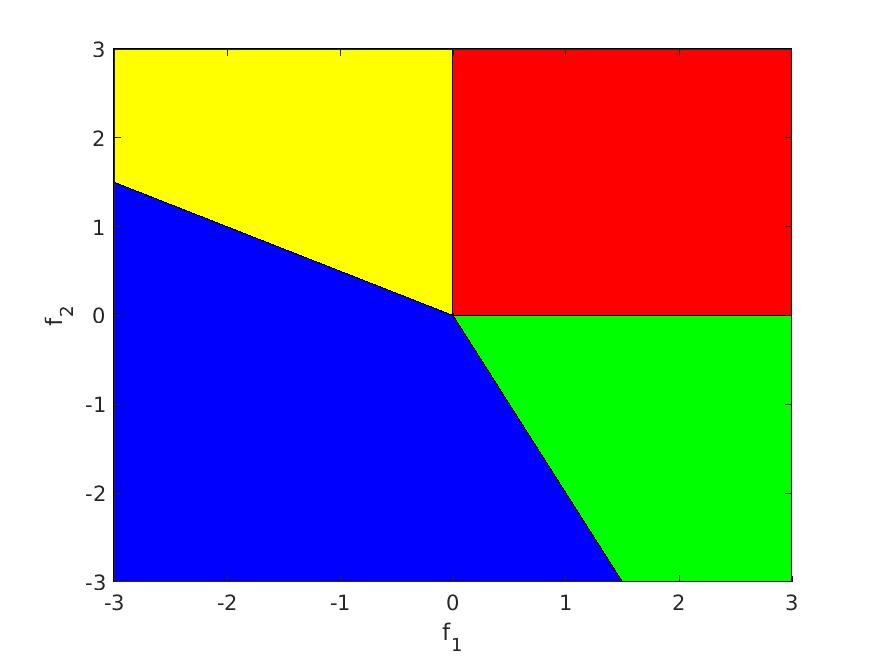}
		\subcaption{\footnotesize Second exactness cone}
	\end{minipage}
	\caption{\label{fig:fourpoint2} Finite set. Left (a): second moment relaxation (light gray) which is the convex hull  of $\sX$ (four black points), and normal cones at the points (colored). Right (b): second exactness cone (colored).}
\end{figure}

On the left of Figure \ref{fig:fourpoint1} we represent the first ($r=1$) moment relaxation $\sR(\bg)^1_1$ (dark gray), the convex hull $\mathrm{conv}\:\sX$ (light gray), and the four points of $\sX$ (black). The tiny dark gray region which remains visible are points in $\sR(\bg)^1_1 \setminus \mathrm{conv}\:\sX$. Also represented (in color) are \new{projections on the coordinates $(f_1,f_2)$ of} the four normal cones at the four points. According to {Theorem} \ref{normal}, up to the sign, they are the four spectrahedral cones $\sS(\hat{\bx})$, $\hat{\bx} \in \sX$ of Theorem \ref{thm}. On the right of Figure \ref{fig:fourpoint1} we represent the exactness cone $\sF(\bg)^1_1$ which is the union of the four spectrahedra, according to Theorem \ref{thm}. We observe tiny conic regions (in white) corresponding to $\R[\bx]_1 \setminus \sF(\bg)^1_1$, namely first degree polynomials $f$ for which the moment-SOS relaxation of first order is not exact. If we solve the relaxation, we hit the slightly inflated tiny regions (dark gray on the left figure) of the moment relaxation $\sR(\bg)^1_1$, yielding a strict lower bound on the value $v(f)$. If instead we minimize the polynomials in $\sF(\bg)^1_1$, we hit one of the four points of $\sX$, i.e. the relaxation is exact.

On Figure \ref{fig:fourpoint2} we represent the same objects for the second relaxation, i.e. $r=2$. On the left, we see that the moment relaxation $\sR(\bg)^2_1$ is the polytope $\mathrm{conv}\:\sX$, i.e. $\sR(\bg)^2_1 \setminus \mathrm{conv}\:\sX$ is empty: the tiny dark gray regions of Figure \ref{fig:fourpoint1}(a) disappeared. We observe on the right that the exactness cone $\sF(\bg)^2_1$ is the whole space $\R[\bx]_1$, i.e. the relaxation is exact everywhere: the tiny white regions of Figure \ref{fig:fourpoint1}(b) disappeared, consistently with Lemma \ref{embed}.
Exactness follows from the property that all non-negative bivariate quartics are SOS. Indeed, the dual problem consists of maximizing $v$ such that $f-v$ is positive on the four points of $\sX$, a linear constraint. But this is equivalent to enforcing that $f-v$ is a degree 4 (i.e. $r=2$) SOS subject to the linear constraint.

\newpage

\subsection{Non-convex set}

	\begin{figure}
	\centering
	\begin{minipage}{0.48\textwidth}
		\includegraphics[width=\textwidth]{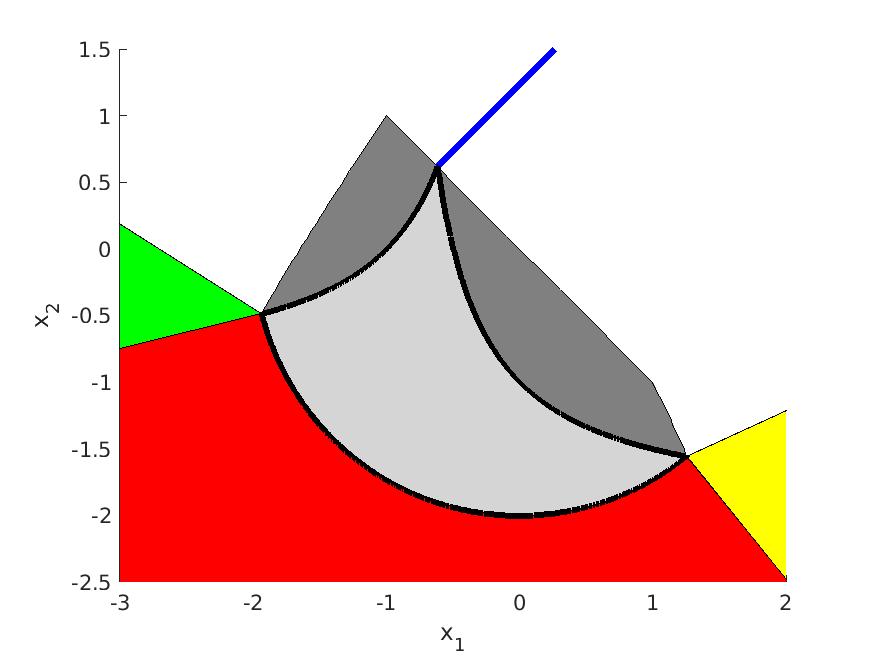}
		\subcaption{\hspace{-.21em}\footnotesize First moment relaxation and  normal cones}
	\end{minipage}
	\begin{minipage}{0.48\textwidth}
		\includegraphics[width=\textwidth]{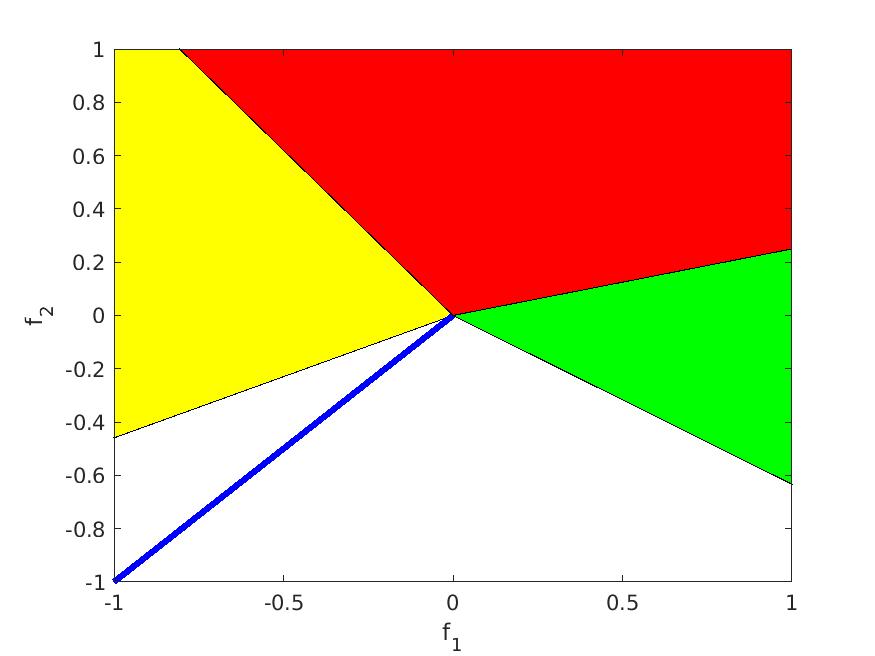}
		\subcaption{\footnotesize First exactness cone}
	\end{minipage}
	\caption{\label{fig:qp1} Non-convex set. Left (a): first moment relaxation (dark gray) of $\sX$ (light gray), and normal cones (colored). Right (b): first exactness cone (colored).}
\end{figure}

\begin{figure}
	\centering
	\begin{minipage}{0.48\textwidth}
		\includegraphics[width=\textwidth]{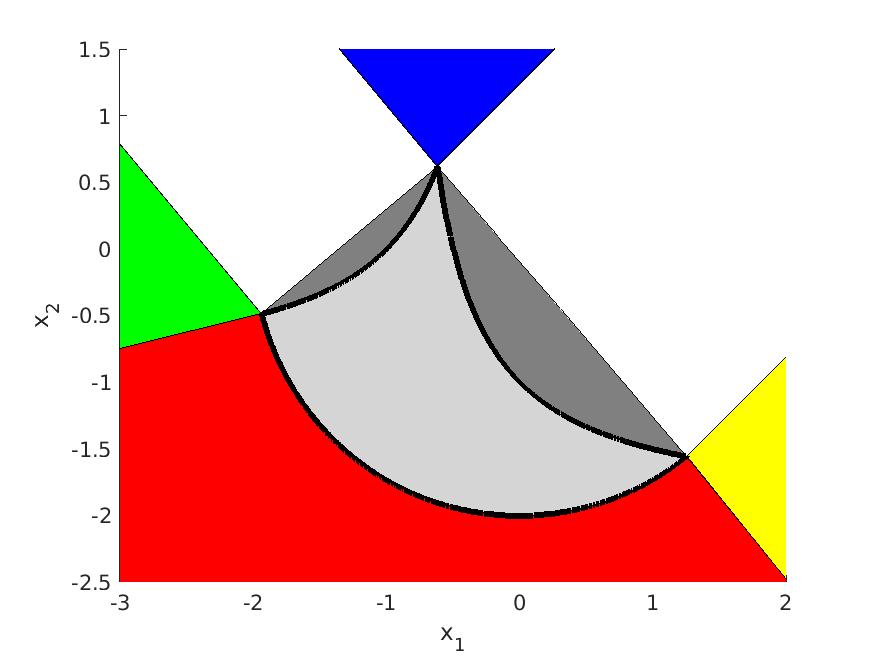}
		\subcaption{\hspace{-.21em}\footnotesize Second moment relaxation and  normal cones}
	\end{minipage}
	\begin{minipage}{0.48\textwidth}
		\includegraphics[width=\textwidth]{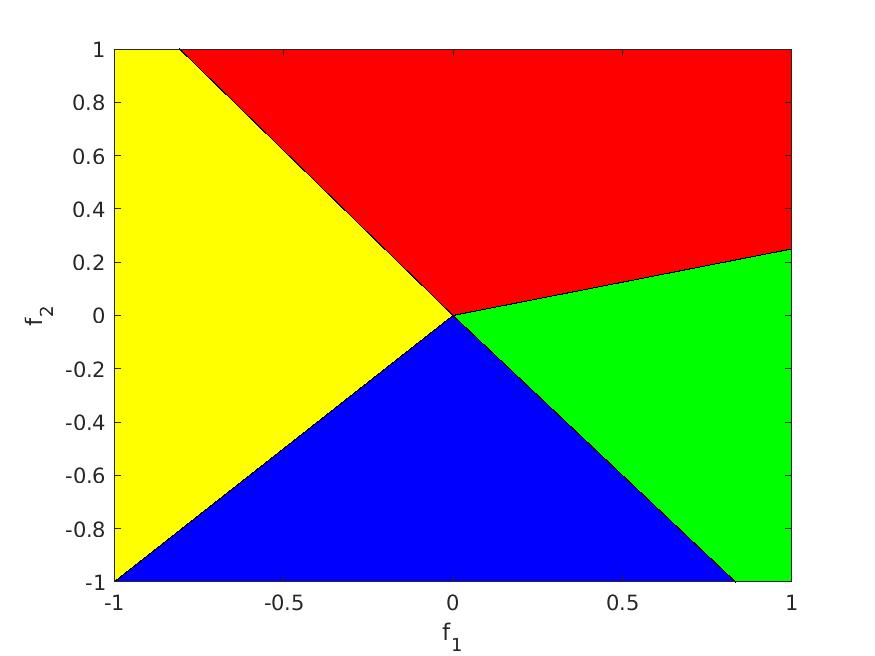}
		\subcaption{\footnotesize Second exactness cone}
	\end{minipage}
	\caption{\label{fig:qp2}  Non-convex set. Left (a): second moment relaxation (dark gray) which is the convex hull  of $\sX$ (light gray), and normal cones at the points (colored). Right (b): second exactness cone (colored).}
\end{figure}

Consider \cite[Ex.2 21]{h23} where $d=1$ and
\[
\begin{array}{rcl}
\sX & = & \{\bx \in \R^2 : g_1(\bx) = 4-x^2_1-x^2_2 \geq 0, \: g_2(\bx)=-1-2x_1-x_2-x_1x_2 \geq 0, \\
& & \quad g_3(\bx)=1+x_1+x_1x_2 \geq 0 \}.
\end{array}
\]
On the left of Figure \ref{fig:qp1} we represent the first ($r=1$) moment relaxation $\sR(\bg)^1_1$ (dark gray) of $\sX$ (light gray), as well as the normal cones to the points of the boundary of $\sR(\bg)^1_1$ where the first relaxation is exact. The green region is the normal cone to the left corner point of $\sX$, the yellow region is the normal cone to the right corner point of $\sX$, the blue line is the one-dimensional normal cone to the top  corner point of $\sX$, and the red region is the union of all the one-dimensional normal cones to the convex circular bottom part of $\sX$. According to {Theorem} \ref{normal}, up to the sign, the green, yellow, and blue cones are spectrahedral cones, whereas the red region is the union of spectrahedral cones along the circular arc. On the right of Figure \ref{fig:qp1} we represent the exactness cone $\sF(\bg)^1_1$ which is the union of these spectrahedral cones, according to Theorem \ref{thm}. The blue line corresponds to the objective function $f(\bx) = -x_1-x_2$ for which the first moment relaxation is exact. It is surrounded by a white region corresponding to objective functions for which the first moment relaxation is not exact. The other colored regions belong to the exactness cone.

On the left of Figure \ref{fig:qp2} we represent the second ($r=2$) moment relaxation $\sR(\bg)^2_1$ (dark gray) of $\sX$ (light gray), as well as the normal cones to the points of the boundary of $\sR(\bg)^2_1$ where the second relaxation is exact. Observe that $\sR(\bg)^2_1 = \mathrm{conv}\:\sX$. In comparison with Figure \ref{fig:qp1}, we notice that the green and yellow normal cones are now larger, and the blue half-line of the first relaxation is now a full-dimensional normal cone to the top corner of $\sX$. On the right of Figure \ref{fig:qp2}, we consistently see that the exactness cone $\sF(\bg)^2_1$ now fills up to whole space $\R[\bx]_1$, i.e. the second moment relaxation is exact for all first degree objective functions.

\section*{Acknowledgement}

The proofs of Theorems \ref{thm} and \ref{normal}, \new{as well as a simplified statement of Definition \ref{exactnesscone},} were suggested by anonymous reviewers. \new{I am also grateful to David de Laat for useful comments.}

\section*{Conflict of interests}

I do not have any conflict of interest to declare regarding this submission.

\end{document}